\documentclass[twoside,12pt,leqno]{amsart}

\usepackage{amsmath}
\usepackage{amssymb}
\usepackage{amscd}

\numberwithin{equation}{section}

\newtheorem{thm}{Theorem}
\newtheorem{lem}[thm]{Lemma}
\newtheorem{prop}[thm]{Proposition}
\newtheorem{cor}[thm]{Corollary}
\newtheorem{remark}[thm]{Remark}

\newcommand\RRR{{\mathbb R}}
\newcommand\FF{{\mathcal F}}
\newcommand\PP{{\mathcal P}}

\newcommand\Int{\mathrm{Int}}
\newcommand\id{\mathrm{id}}

\newcommand\Orbit{\mathcal{O}}
\newcommand\Stab{\mathcal{S}}
\newcommand\Diff{\mathcal{D}}

\newcommand\StabId{\Stab_{\id}}
\newcommand\Stabfcr{\Stab(\mrsfunc,\singf)}
\newcommand\StabIdfcr{\StabId(\mrsfunc,\singf)}
\newcommand\Stabf{\Stab(\mrsfunc)}
\newcommand\StabIdf{\Stab_{\id}(\mrsfunc)}

\newcommand\DiffM{\Diff(\manif)}
\newcommand\DiffIdM{\Diff_{\id}(\manif)}
\newcommand\DiffMcr{\Diff(\mrsfunc,\singf)}
\newcommand\DiffIdMcr{\Diff_{\id}(\mrsfunc,\singf)}

\newcommand\Orbf{\Orbit(\mrsfunc)}
\newcommand\Orbff{\Orbit_{\mrsfunc}(\mrsfunc)}
\newcommand\Orbfcr{\Orbit(\mrsfunc,\singf)}
\newcommand\Orbffcr{\Orbit_{\mrsfunc}(\mrsfunc,\singf)}

\newcommand\Psp{P}
\newcommand\manif{M}
\newcommand\mrsfunc{f}
\newcommand\gfunc{g}
\newcommand\difM{h}
\newcommand\singf{\Sigma_\mrsfunc}

\newcommand\Korbf{\FF(\mrsfunc)}
\newcommand\ConfSpnM{\FF_{n}(\Int\manif)}
\newcommand\cnfpr{e}
\newcommand\PermutGr[1]{\mathbb{S}_{#1}}
\newcommand\DiffFiber{\hat\Diff(\mrsfunc)}
\newcommand\OrbFiber{\hat\Orbit(\mrsfunc)}
\newcommand\gdim{{\mathrm{g.d.}}}
\newcommand\hdim{{\mathrm{h.d.}}}

\author{Sergiy Maksymenko}
\title{Homotopy dimension of orbits of Morse functions on surfaces}

\begin{document}
\maketitle

\begin{abstract}
Let $\manif$ be a compact surface, $\Psp$ be either the real line $\RRR$ or the circle $S^1$, and $\mrsfunc:\manif\to\Psp$ be a $C^{\infty}$ Morse map.
The identity component $\DiffIdM$ of the group of diffeomorphisms of $\manif$ acts on the space $C^{\infty}(\manif,\Psp)$ by the following formula:
$h\cdot \mrsfunc = \mrsfunc \circ \difM^{-1}$ for $\difM\in\DiffIdM$ and $\mrsfunc\in C^{\infty}(\manif,\Psp)$.
Let $\Orbf$ be the orbit of $\mrsfunc$ with respect to this action and $n$ be the total number of critical points of $\mrsfunc$.
In this note we show that $\Orbf$ is homotopy equivalent to a certain covering space of the $n$-th configuration space of the interior $\Int\manif$.
This in particular implies that the (co-)homology of $\Orbf$ vanish in dimensions greater than $2n-1$, and the fundamental group $\pi_1\Orbf$ is a subgroup of the $n$-th braid group $\mathcal{B}_{n}(\manif)$.
\end{abstract}


{\bf AMS Classification 2000:} 
14F35, % Homotopy theory; fundamental groups 
46T10 % Manifolds of mappings 

%%%%%%%%%%%%%%%%%%%%%%%%%%%%%%%%%
\section{Introduction}
Let $\manif$ be a compact surface, $\Psp$ be either the real line $\RRR$ or the circle $S^1$.
Then the group $\DiffM$ of $C^{\infty}$ diffeomorphisms of $\manif$ acts on the space $C^{\infty}(\manif,\Psp)$ by the following formula:
\begin{equation}\label{equ:act-DM}
h\cdot \mrsfunc = \mrsfunc \circ \difM^{-1}
\end{equation}
for $\difM\in\DiffM$ and $\mrsfunc\in C^{\infty}(\manif,\Psp)$.

We say that a smooth ($C^{\infty}$) map $\mrsfunc:\manif\to\Psp$ is \emph{Morse\/} if 
\begin{enumerate}
 \item[(i)] critical points of $\mrsfunc$ are non-degenerate and belong to the interior of $\manif$;
 \item[(ii)] $\mrsfunc$ is constant on every connected component of $\partial\manif$.
\end{enumerate}

Let $\mrsfunc\in C^{\infty}(\manif,\Psp)$, $\singf$ be the set of critical points of $\mrsfunc$, and $\DiffMcr$ be the subgroup of $\DiffM$ consisting of diffeomorphisms $\difM$ such that $\difM(\singf)=(\singf)$.

Then we can define the stabilizers $\Stabf$ and $\Stabfcr$, and orbits $\Orbf$ and $\Orbfcr$ with  respect to the actions of the groups $\DiffM$ and $\DiffMcr$.
Thus 
$$
\Stabf=\{\difM\in\DiffM \ : \ \mrsfunc\circ\difM^{-1}=\mrsfunc  \},
\qquad 
\Orbf=\{\mrsfunc \circ \difM^{-1} \ : \ \difM\in\DiffM  \},
$$
$$
\Stabfcr = \Stabf \cap \DiffMcr.
$$

We endow the spaces $\DiffM$ and $C^{\infty}(\manif,\Psp)$ with the corresponding $C^{\infty}$ Whitney topologies.
They induce certain topologies on the stabilizers and orbits.

Let $\DiffIdM$ and $\DiffIdMcr$ be the identity path components of the groups $\DiffM$ and $\DiffMcr$, 
$\StabIdf$ and $\StabIdfcr$ be the identity path components of the corresponding stabilizers, and $\Orbff$ and $\Orbffcr$ be the path-components of $\mrsfunc$ in the corresponding orbits with respect to the induced topologies.

\begin{lem}
If $\singf$ is discrete set, e.g. when $\mrsfunc$ is Morse, then $\StabIdfcr=\StabIdf$.
\end{lem}
\begin{proof}
Since $\Stabfcr \subset \Stabf$, we have that $\StabIdfcr\subset\StabIdf$.
Conversely, let $h_t:\manif\to\manif$ be an isotopy such that $h_0=\id_{\manif}$ and $h_t\in\Stabf$ for all $t\in I$., i.e. $\mrsfunc\circ h_t=\mrsfunc$.
We have to show that $h_t\in\Stabfcr$ for all $t\in I$.
Notice that $d(\mrsfunc\circ h_t)=h_t^{*}d\mrsfunc=d\mrsfunc$, whence $h_t(\singf)=\singf$.
Since $\singf$ is discrete and $h_0=\id_{\manif}$ fixes $\singf$, we see that so does every $h_t$, i.e. $h_t\in\Stabfcr$.
\end{proof}

Let $\mrsfunc:\manif\to\Psp$ be a Morse map.
Denote by $c_i$, $(i=0,1,2)$, the total numbers of critical points of $\mrsfunc$ of index $i$ and let $n=c_0+c_1+c_2$ be the total number of critical points of $\mrsfunc$.

Notice that for every Morse map $\mrsfunc$  its orbits $\Orbf$ and $\Orbfcr$ are Fr\'echet submanifolds of $C^{\infty}(\manif,\Psp)$ of finite codimension, see~\cite{Poenaru, Sergeraert}.
Therefore, e.g.~\cite{Palais}, these orbits have the homotopy types of CW-complexes.
But in general these complexes may have infinite dimensions.

Let $X$ be a topological space which is homotopy equivalent to some CW-complex.
Then a \emph{homotopy dimension} $\hdim\,X$ of $X$ is the minimal dimension of a CW-complex homotopy equivalent to $X$.
In particular $\hdim\,X$ can be equal to $\infty$.
It is also evident that if $\hdim\,X<\infty$, then (co-)homology of $X$ vanish in dimensions greater that $\hdim\,X$.

If $\pi$ is  a finitely presented group $\pi$, then the \emph{geometric dimension\/} of $\pi$, denoted $\gdim\,\pi$, is the homotopy dimension of its Eilenberg-Mac Lane space $K(\pi,1)$:
$$
\gdim\,\pi := \hdim\, K(\pi,1).
$$

In~\cite[Theorems\;1.3,\;1.5,\;1.9]{Maks:AGAG:2006} the author described the homotopy types of $\StabIdf$, $\Orbff$, and $\Orbffcr$.
It follows from these results that
$$\hdim\,\StabIdf, \;\, \; \hdim\,\Orbffcr\;\;\leq\; \; 1.$$
In fact, $\StabIdf$ is contractible provided either $\mrsfunc$ has at least one critical point of index $1$, i.e., $c_1\geq1$ or $\manif$ is non-orientable.
Otherwise $\StabIdf \simeq S^1$.

Also, $\Orbff\simeq S^1$ for Morse mappings $T^2\to S^1$ and $K^2 \to S^1$ without critical points, and $\Orbff$ is contractible in all other cases, where $K$ stands for the Klein bottle.

For $\Orbff$ the description is not so complete.
But if $\mrsfunc$ is \emph{generic},  i.e., it takes distinct values at distinct critical points, then 
$$\hdim\,\Orbff \;\; \leq \;\; \max\{c_0 + c_2 + 1, \,  c_1+2\} \;\;<\;\; \infty.$$
Actually, in this case $\Orbff$ is either contractible or homotopy equivalent to $T^k$ or to $\RRR{P}^3\times T^k$ for some $k\geq0$, where $T^k$ is a $k$-dimensional torus.

Thus the upper bound for $\hdim\,\Orbff$ (at least in generic case) depends only on the number of critical points of $\mrsfunc$ at each index.

In this note we will show that $\hdim\,\Orbff \leq 2n-1$ for arbitrary Morse mapping $\mrsfunc:\manif\to\Psp$ having exactly $n\geq1$ critical points.
Notice that if $n=0$, then $\mrsfunc$ is generic, and in fact $\hdim\,\Orbff \leq 1$, see~\cite[Table\;1.10]{Maks:AGAG:2006}.

\begin{thm}\label{th:hom-type-Orbits}
Let $\mrsfunc:\manif\to\Psp$ be a Morse map and $n$ be the total number of critical points of $\mrsfunc$.
Assume that $n\geq1$.
Denote by $\ConfSpnM$ the configuration space of $n$ points of the interior $\Int\manif$ of $\manif$.
Then $\Orbff$ is homotopy equivalent to a certain covering space $\Korbf$ of $\ConfSpnM$.
\end{thm}

\begin{cor}
$\hdim\Orbff \leq 2n-1$, whence (co-)homology of $\Orbff$ vanish in dimensions $\geq 2n$.
\end{cor}
\begin{proof}
Since $\ConfSpnM$ and its connected covering spaces are \emph{open\/} manifolds of dimension $2n$, they are homotopy equivalent to CW-complexes of dimensions not greater than $2n-1$.
\end{proof}

For simplicity denote $\pi = \pi_1\Orbff$.
Since the covering map $\Korbf \to \ConfSpnM$ yields a monomorphisms of fundamental groups, we obtain the following: 

\begin{cor}
The fundamental group $\pi$ of $\Orbff$ is a subgroup of the $n$-th braid group $\mathcal{B}_{n}(\manif) = \pi_1(\ConfSpnM)$ of $\manif$.
\end{cor}
\begin{cor}\label{cor:asph-orb}
%Suppose that $\mrsfunc$ has at least one critical point of index $1$ and $\manif$ is aspherical, i.e., 
Suppose that $\manif$ is aspherical, i.e., $\manif\not= S^2, \RRR P^2$.
Then $\Orbff$ is aspherical as well, i.e., $K(\pi,1)$-space, whence \ $\gdim\, \pi \leq 2n-1.$
\end{cor}
\begin{proof}
Actually the aspherity of $\Orbff$ for the case $\manif\not= S^2, \RRR P^2$ is proved in~\cite[Theorems~1.5, 1.9]{Maks:AGAG:2006}.

But it can be shown by another arguments.
It is well known and can easily be deduced from~\cite{FadellNeuwirth}  that for an aspherical surface $\manif$ every of its configuration spaces $\ConfSpnM$ and thus every covering space of $\ConfSpnM$  are aspherical as well.
Hence so is $\Korbf$ and thus $\Orbff$ itself.
\end{proof}

A presentation for $\pi$ will be given in another paper.

\section{Orbits of the actions of $\DiffIdM$ and $\DiffIdMcr$}

\begin{prop}\label{pr:Orbff_is_orbit_for_DiffIdM}
Let $\mrsfunc:\manif\to\Psp$ be a Morse map and 
\begin{equation}\label{equ:p-d2o}
p:\DiffM\mapsto\Orbf, \qquad p(\difM) =\mrsfunc\circ\difM^{-1}
\end{equation}
be the natural projection.
Then $\Orbff$ is the orbit of $\mrsfunc$ with respect to $\DiffIdM$ and $\Orbffcr$ is the orbit of $\mrsfunc$ with respect to $\DiffIdMcr$.
In other words,
$$ p(\DiffIdM)=\Orbff \qquad\text{and}\qquad  p(\DiffIdMcr)=\Orbffcr. $$
\end{prop}
\begin{proof}
The proof is based on the following general statement.
Let $G$ be a topological group transitively acting on a topological space $O$ and $f\in O$.
Denote by $G_e$ the \emph{path-component\/} of the unit $e$ in $G$ and let $O_f$ be the \emph{path-component\/} of $f$ in $O$.
\begin{lem}\label{lm:p_Ge_Xf}
Suppose that the mapping $p:G\to O$ defined by 
$$p(\gamma)=\gamma\cdot f, \qquad  \forall \gamma\in G$$
 satisfies a covering path axiom (in particular, this holds when $p$ is a locally trivial fibration).
Then $O_f$ is the orbit of $f$ with respect to the induced action of $G_e$ on $O$, i.e., $p(G_e)=O_f$.
\end{lem}
\proof
Evidently, $p(G_e)\subset O_f$.
Conversely, let $g\in O_f$. 
Then there exists a path $\omega:I\to O_f$ between $f$ and $g$, i.e., $\omega(0)=f$ and $\omega(1)=g$.
Since $p$ satisfies the covering path axiom, $\omega$ lifts to the path $\tilde\omega:I\to G$ such that $\tilde\omega(0)=e$ and $\omega=p\circ \tilde\omega$.
Then $g= \omega(1) = p\circ \tilde\omega(1) \in p(G_e)$.
Thus $p(G_e)=O_f$.
\endproof
It remains to note that the mapping~\eqref{equ:p-d2o} is a locally trivial fibration, see e.g.~\cite{Poenaru,Sergeraert}, and $\DiffM$ (resp. $\DiffMcr$) transitively acts on the orbit $\Orbf$ (resp. $\Orbfcr$).
Therefore the conditions of Lemma~\ref{lm:p_Ge_Xf} are satisfied.
\end{proof}

\section{Proof of Theorem~\ref{th:hom-type-Orbits}}
Let $\ConfSpnM$ be the configuration space of $n$ points of the interior $\Int\manif$ of $M$.
Thus 
\begin{equation}\label{equ:conf_sp}
\ConfSpnM = \PP_{n}(\Int\manif)/\PermutGr{n}, 
\end{equation}
where 
$$
\PP_{n}(\Int\manif) = \{ (x_1,\ldots x_n) \ | \ x_i\in \Int\manif \ \text{and} \ x_i\not=x_j \ \text{for} \ i\not=j  \}
$$
is called the \emph{pure} $n$-th configuration space of $\Int\manif$, and $\PermutGr{n}$ is the symmetric group of $n$ symbols freely acting on $\PP_{n}(\Int\manif)$ by permutations of coordinates.

We can regard $\ConfSpnM$ as the space of $n$-tuples of mutually distinct points of $\Int\manif$.

Denote by $\singf=\{x_1,\ldots,x_n\}$ the set of critical points of $\mrsfunc$. 
Then for every $\gfunc\in\Orbff$ the set $\Sigma_{\gfunc}$ of its critical points is a point in $\ConfSpnM$.
Hence the correspondence $\gfunc\mapsto\Sigma_{\gfunc}$ is a well-defined mapping 
$$
k:\Orbff\to\ConfSpnM, \qquad 
k(\gfunc)=\Sigma_{\gfunc}.
$$

\begin{lem}\label{lm:prop_of_k}
{\rm(i)}~The mapping $k$ is a locally trivial fibration.
The connected component of the fiber containing $\mrsfunc$ is homeomorphic to $\Orbffcr$.

{\rm(ii)}~Let $k_{i}:\pi_{i}(\Orbff,f)\to\pi_{i}(\ConfSpnM,\singf)$, $(i\geq1)$, be the corresponding homomorphism of homotopy groups induced by $k$.
Then $k_1$ is a {\bfseries monomorphism} and all other $k_i$ for $i\geq2$ are isomorphisms.
\end{lem}
Assuming that Lemma~\ref{lm:prop_of_k} is proved we will now complete our theorem.
Let $\Korbf$ be the covering space of $\ConfSpnM$ corresponding to the subgroup 
$$\pi_1\Orbff \approx k_1(\pi_1\Orbff) \subset \pi_1 \ConfSpnM.$$

Then $k$ lifts to the mapping $\hat k:\Orbff\to\Korbf$ which induces isomorphism of all homotopy groups.
Since $\Orbff$ and $\Korbf$ are connected, we obtain from (2) that $\hat k$ is a desired homotopy equivalence.
Theorem~\ref{th:hom-type-Orbits} is proved modulo Lemma~\ref{lm:prop_of_k}.

\begin{proof}[Proof of Lemma~\ref{lm:prop_of_k}.]
(i) Recall, \cite{FadellNeuwirth}, that the following \emph{evaluation\/} map 
$$\cnfpr:\DiffIdM\to\ConfSpnM,\qquad \cnfpr(\difM)=\difM(\singf)$$
is a locally trivial principal fibration with fiber $$\DiffFiber=\DiffIdM\cap\DiffMcr.$$

Let $p:\DiffIdM\to\Orbff$ be the projection defined by $p(\difM)=\mrsfunc\circ\difM^{-1}$.
Then the set of critical points of the function $\mrsfunc\circ\difM^{-1} \in\Orbff$ is $\difM(\singf)$.
Therefore $e$ coincides with the following composition:
$$
\begin{CD}
\cnfpr = k \circ p : 
\DiffIdM @>{p}>> \Orbff @>{k}>> \ConfSpnM.
\end{CD}
$$
Since $\cnfpr$ and (by Proposition~\ref{pr:Orbff_is_orbit_for_DiffIdM}) the mapping $p$ are principal locally trivial fibrations, we obtain that $k$ is also a locally trivial fibration with fiber $\OrbFiber$ being the orbit of $\mrsfunc$ with respect to the group $\DiffFiber$.

It is easy to see that the identity component of the group $\DiffFiber$ coincides with $\DiffIdMcr$, whence by Proposition~\ref{pr:Orbff_is_orbit_for_DiffIdM}, the connected component of $\OrbFiber$ containing $\mrsfunc$ is $\Orbffcr$.

(ii)
As noted above since $n\geq1$, it follows from~\cite[Theorems~1.5(i), 1.9]{Maks:AGAG:2006} that $\Orbffcr$ is contractible.
Then from the exact sequence of homotopy groups of the fibration $k$ we obtain that for $i\geq 2$ every $k_i$ is an isomorphism, and $k_1$ is a monomorphism.
Lemma~\ref{lm:prop_of_k} is proved.
\end{proof}

\begin{remark}\rm
In general the covering map $\Korbf \to \ConfSpnM$ is not \emph{regular}, i.e., $\pi_1\Orbff\approx\pi_1\Korbf$ is not a normal subgroup of $\mathcal{B}_{n}(\manif)= \pi_1\ConfSpnM$.
\end{remark}

\begin{remark}\rm
Theorem~\ref{th:hom-type-Orbits} does not answer the question whether $\Orbff$ has the homotopy type of a \emph{finite\/} CW-complex.
Indeed, since $\manif$ is compact, it follows from~\eqref{equ:conf_sp} that $\mathcal{B}_{n}(\manif)$ can be regarded as an open cellular (i.e.\! consisting of full cells) subset of a finite CW-complex $\prod_{n}\manif/\PermutGr{n}$.
Therefore if the covering map $\Korbf\to\ConfSpnM$ is an infinite sheet covering, i.e., $\pi_1\Orbff$ has an infinite index in $\mathcal{B}_{n}(\manif)$, then we obtain a priori an infinite cellular subdivision of $\Korbf$.
On the other hand, as noted above, for a generic Morse map $\mrsfunc:\manif\to\Psp$ a finiteness of the homotopy type of $\Orbff$ follows from~\cite{Maks:AGAG:2006}.
\end{remark}

%%%%%%%%%%%%%%%%%%%%%%%%%%%%%%%%%%%%%%%%%%%%%%%%%%%%%%%%%%%%%%%%%%%%%%%
%%%%%%%%%%%%%%%%%%%%%%%%%%%%%%%%%%%%%%%%%%%%%%%%%%%%%%%%%%%%%%%%%%%%
%%%%%          The bibliography
%%%%%%%%%%%%%%%%%%%%%%%%%%%%%%%%%%%%%%%%%%%%%%%%%%%%%%%%%%%%%%%%%%%%%

%%%%%%%%%%%%%%%%%%%%%%%%%%%%%%%%%%%%%%%%%%%%%%%%%%%%%%%%%%%%
%%%%%%%%%%%%%%%%%%%%%%%%%%%%%%%%%%%%%%%%%%%
%% The address
\noindent 
Sergiy Maksymenko \\
Topology dept., Institute of Mathematics of NAS of Ukraine, \\ 
Te\-re\-shchenkivska st. 3, Kyiv, 01601 Ukraine \\
E-mail: \texttt{maks@imath.kiev.ua}\\
%%%%%%%%%%%%%%%%%%%%%%%%%%%%%%%%%%%%%%%%%%%%%%%%%%%%%%%%%%%%%%%%%

\end{document}